\documentclass[10pt]{article} 

\usepackage[titles]{tocloft} 
\usepackage{setspace}
\usepackage{float,url}
\usepackage[utf8]{inputenc} 

\usepackage[margin=1in]{geometry}
\title{A Stochastic Programming Approach to the Railcar Maintenance Problem with Service Level and Track Capacity Considerations}

\author{Murat Elh\"{u}seyni\thanks{murat.elhuseyni@yeditepe.edu.tr, Industrial Engineering, Yeditepe University, 34755, Istanbul, Turkey.  Corresponding author.
} \and
Burak Kocuk\thanks{burakkocuk@sabanciuniv.edu, Industrial Engineering,  Sabanc{\i} University, 34956,  Istanbul, Turkey.
}}

\usepackage{tikz}
\usepackage{pgfplots}
\usepackage{verbatim}
\usepackage{bigstrut}

\usepackage{algorithm}
\usepackage[noend]{algorithmic}

\usepackage{subfig} 
\usepackage{caption}
\pgfplotsset{compat=1.18}
\usepackage{graphicx}
\usepackage{epstopdf}

\usepackage{booktabs}
\usepackage{array} 
\usepackage{paralist} 
\usepackage{verbatim} 

\usepackage{fancyhdr} 
\usepackage{sectsty}
\allsectionsfont{\sffamily\mdseries\upshape} 

\usepackage[nottoc,notlof,notlot]{tocbibind} 
\usepackage[titles]{tocloft} 

\usepackage{amsmath, amsthm, amssymb}

\usepackage{longtable}
\usepackage{pdflscape}
\usepackage{natbib}

\usepackage{adjustbox}
\usepackage{multirow}

\usepackage[margin=1in]{geometry}
\usepackage{amssymb}
\usepackage{float}
\usepackage{color}
\usepackage{soul}
\usepackage{mathtools}
\let\hat\widehat

\usepackage{enumitem}

\usepackage{amsthm}

\allowdisplaybreaks

\DeclareMathOperator*{\argmin}{argmin}

\usepackage{rotating}
\usepackage{pdflscape}

\usepackage{pgfplots}
\usepackage{pgfplotstable} 
\usepgfplotslibrary{groupplots}
\usepgfplotslibrary{units}
\usepackage{url}
\usepackage{color}

\graphicspath{{figures/}} 

\allowdisplaybreaks

\usepackage{color}
\usepackage{bm}
\usepackage{bbm}

  \pgfplotsset{compat=1.18} 

\usepackage{pgfgantt}

\begin{document}

\date{} 
\maketitle

 \abstract{
 Railcars, as part of the rolling stock, perform regular transportation tasks with respect to a service level agreement (SLA) and undergo preventive maintenance at regular intervals based on the recommendations of train manufacturers. When unexpected failures occur, they need to enter corrective maintenance immediately. However, this reactive approach may result in large SLA violations and an excessive number of corrective maintenance actions. In this study, we utilize a predictive maintenance approach based on the reliability of a railcar.  In particular, we propose a stochastic programming model, in which railcar failure scenarios are generated from a Weibull distribution, a common assumption in the reliability literature. The model incorporates both SLA and track-capacity considerations and is solved through the Sample Average Approximation (SAA) method. We generate random instances to compare the stochastic model and a deterministic model adopted from the literature with respect to several system parameters. Our results show that the stochastic model achieves lower total costs, fewer SLA violations, and a reduced number of corrective interventions compared with deterministic approaches, while effectively managing track-capacity constraints. Our results underscore the importance of the predictive approach in the context of the railcar maintenance problem.
}





\section{Introduction} \label{s:introduction}

Rapid urbanization in the 21st century has resulted in a sharp rise in transportation demand across metropolitan areas. To meet this demand, urban rail service providers are under severe pressure to operate efficiently. A holistic approach that carefully manages both operational and  maintenance activities is essential to ensure uninterrupted train service. However, in practice, operators often prioritize immediate operational goals at the expense of broader maintenance strategies. This imbalance may result in deferred maintenance tasks, elevated corrective intervention rates, and unexpected service disruptions. Over time, the accumulation of such short-term decisions escalates major repair costs, undermines system reliability, and threatens transportation safety. These effects not only make it increasingly difficult to sustain high service quality but also erode the reputation of service providers. Addressing these challenges necessitates structured frameworks that integrate service operations with maintenance planning, thereby enabling operators to balance immediate performance requirements against long-term asset health. Within this context, a central question arises: how should operators allocate maintenance resources between preventive and corrective actions while still fulfilling service obligations? To explore this issue, we examine a rail service system managed by an operator who must simultaneously address operational and maintenance decisions.

With regard to the operational aspect, the system operator is bound by a \textit{Service Level Agreement} (SLA), which specifies the contractual service expectations between the provider and the customer \cite{Cips2009}. In the context of our study, SLA refers to daily passenger service requirements that must be satisfied. On the maintenance side, the operator can perform both preventive and corrective actions: preventive maintenance refers to proactively scheduled interventions aimed at avoiding failures, whereas corrective maintenance refers to reactive interventions carried out after unexpected railcar breakdowns. These two approaches are inherently in conflict. Relying solely on corrective actions exposes the system to unplanned outages and long disruptions, since corrective maintenance generally requires more time than preventive measures. Conversely, scheduling excessive preventive maintenance reduces fleet availability and undermines system efficiency. Hence, the operator must carefully determine the timing and frequency of preventive actions to manage the transportation system cost-effectively while satisfying service obligations (i.e., SLA requirements) under capacity constraints (e.g., track availability).

The standard approach in the reliability literature assigns each railcar a \textit{maintenance interval}, defined by an earliest and latest allowable start time, which are typically derived from manufacturer recommendations (see, e.g.,~\cite{elhuseyni2021integrated,folco2024rolling,liao5026899scenario}). A critical deadline, referred as the \textit{due time}, marks the latest point by which maintenance must begin to avoid increased failure risk. Breaching this interval significantly raises the probability of railcar breakdowns. However, in practice, service providers such as MetroIstanbul, {the primary subway operator in Istanbul}, 
may still deploy railcars beyond this point if doing so is necessary to fulfill SLA requirements, thereby accepting the risk of failure to sustain service levels~\cite{elhuseyni2021integrated}. 

Instead of this standard preventive maintenance approach, we propose a predictive approach in this paper. In particular, assuming that historical failure data is available, we estimate remaining useful life (RUL) distributions and generate plausible failure scenarios. 
This enables the development of a predictive maintenance strategy: railcars can be withdrawn for maintenance based on their reliability status before failures occur. 


Based on the system components defined above, the main goal in this study is to obtain a joint operational and maintenance schedule over a given planning horizon. The objective is to minimize a composite cost function that includes preventive and corrective maintenance costs, penalties associated with SLA and track capacity violations, and operational expenses. The optimization is carried out using inputs such as SLA requirements, maintenance depot capacities, railcar age profiles, and failure scenarios derived from reliability distributions.

\begin{figure}[H]
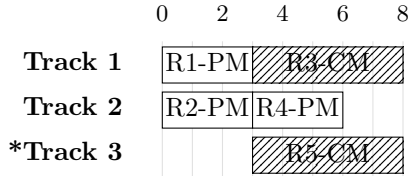

\centering
\begin{ganttchart}[
		x unit=0.4cm,
		y unit chart=0.6cm,
		canvas/.style={draw=none,fill=none}, 
		vgrid={*1{draw=black!12}},           
		inline,                              
		group/.style={draw=none,fill=none},  
		bar top shift=0.125,                   
		bar height=0.8,                      
		y unit title=0.5cm,                  
		title/.style={draw=none,fill=none},  
		include title in canvas=false        
	]{-1}{8}

	\gantttitle{0}{2}
	\gantttitle{2}{2}
	\gantttitle{4}{2}
	\gantttitle{6}{2}
	\gantttitle{8}{2} \\

	\ganttgroup[inline=false]{Track 1}{0}{1}
	\ganttbar[bar/.style={fill=white,draw}]{R1-PM}{0}{2}
	\ganttbar[bar/.style={fill=white,draw, pattern=north east lines}]{R3-CM}{3}{7} \\

	\ganttgroup[inline=false]{Track 2}{0}{1}
	\ganttbar[bar/.style={fill=white,draw}]{R2-PM}{0}{2}
	\ganttbar[bar/.style={fill=white,draw}]{R4-PM}{3}{5} \\

	\ganttgroup[inline=false]{*Track 3}{0}{1}
	\ganttbar[bar/.style={fill=white,draw, pattern=north east lines}]{R5-CM}{3}{7}

\end{ganttchart}
\caption{A maintenance schedule example with two tracks (the horizontal axis represents the time). Track 3 represents the extra capacity. The terms ``Railcar", ``Preventive Maintenance" and ``Corrective Maintenance" are abbreviated as ``R", ``PM" and ``CM", respectively.}
\label{fig:ex-maintenance}
\end{figure}
\begin{figure}[H]
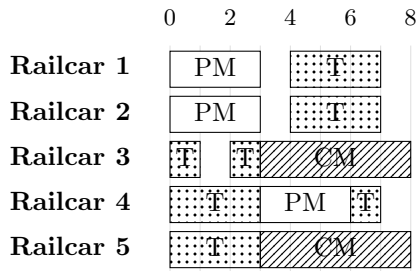

\centering
\begin{ganttchart}[
		x unit=0.4cm,
		y unit chart=0.6cm,
		canvas/.style={draw=none,fill=none}, 
		vgrid={*1{draw=black!12}},           
		inline,                              
		group/.style={draw=none,fill=none},  
		bar top shift=0.125,                   
		bar height=0.8,                      
		y unit title=0.5cm,                  
		title/.style={draw=none,fill=none},  
		include title in canvas=false        
	]{-1}{8}

	\gantttitle{0}{2}
	\gantttitle{2}{2}
	\gantttitle{4}{2}
	\gantttitle{6}{2} 
	\gantttitle{8}{2} \\

	\ganttgroup[inline=false]{Railcar 1}{0}{1}
	\ganttbar[bar/.style={fill=white,draw}]{PM}{0}{2}
	\ganttbar[bar/.style={fill=white,draw, pattern=dots}]{T}{4}{6} 
 \\

	\ganttgroup[inline=false]{Railcar 2}{0}{1}
	\ganttbar[bar/.style={fill=white,draw}]{PM}{0}{2}
	\ganttbar[bar/.style={fill=white,draw, pattern=dots}]{T}{4}{6} 
 \\

	\ganttgroup[inline=false]{Railcar 3}{0}{1}
 	\ganttbar[bar/.style={fill=white,draw, pattern=dots}]{T}{0}{0}
   	\ganttbar[bar/.style={fill=white,draw, pattern=dots}]{T}{2}{2}
	\ganttbar[bar/.style={fill=white,draw, pattern=north east lines}]{CM}{3}{7}
 \\

	\ganttgroup[inline=false]{Railcar 4}{0}{1}
	\ganttbar[bar/.style={fill=white,draw, pattern=dots}]{T}{0}{2}
	\ganttbar[bar/.style={fill=white,draw}]{PM}{3}{5}
	\ganttbar[bar/.style={fill=white,draw, pattern=dots}]{T}{6}{6}
 \\

	\ganttgroup[inline=false]{Railcar 5}{0}{1}
	\ganttbar[bar/.style={fill=white,draw, pattern=dots}]{T}{0}{2}
	\ganttbar[bar/.style={fill=white,draw, pattern=north east lines}]{CM}{3}{7}

\end{ganttchart}
\caption{A railcar utilization schedule example  (the horizontal axis represents the time).  The term ``Transport Task" is abbreviated as~``T".}
\label{fig:ex-utilization}
\end{figure}

To illustrate the proposed scheduling framework, we provide two figures—Figure~\ref{fig:ex-maintenance} and Figure~\ref{fig:ex-utilization}—which depict the maintenance and utilization plans for a scenario involving five railcars and two tracks over a two-day period. Each day is divided into four periods, and the SLA follows a recurring demand pattern of 3–2–3–0. The first and third periods correspond to peak demand, the second to moderate demand, and the fourth (night) to no transport service. In this example, the depot normally operates with two available maintenance tracks. {The railcars are assumed to experience breakdowns at the following times, listed in ascending order of railcar index: 3, 2, 3, 4, and 3.} To mitigate failures, maintenance is scheduled proactively for railcars 1, 2, and 4, all of which receive preventive maintenance before their respective breakdown times. In contrast, railcars 3 and 5 are scheduled for maintenance at or after their failure times. As a result, these railcars undergo immediate corrective maintenance at time 3. To accommodate this reactive maintenance, an additional (third) track is temporarily opened, incurring a capacity expansion cost.


From the perspective of the operations department, Figure~\ref{fig:ex-utilization} presents a period-by-period Gantt chart of each railcar’s activity: transport (T), preventive maintenance (PM), corrective maintenance (CM), or idle time. For example, Railcar~3 remains idle during period~2. Notably, during period~5, the SLA requires three railcars in operation, but only two are available. This results in one SLA violation, which incurs an associated penalty in the cost function.

%





\subsection{Literature Review}

Research on railcar preventive maintenance within train timetables often centers around the challenge of directing railcars to maintenance facilities. One of the earliest contributions in this line of work is the study by Maroti and Kroon \cite{maroti2007maintenance}, which proposes reassigning urgent train units to appropriate maintenance depots by swapping them with non-urgent units in the existing timetable, thereby ensuring timely maintenance without disrupting service continuity. In related work on train unit circulation, Giacco et al. \cite{Giacco2014b} focuses on minimizing the number of train units and maintenance tasks, while evaluating whether empty train services should be limited. Lai et al. \cite{Lai2015} address the scheduling of daily and monthly inspections, and later extended this work to high-speed trains in \cite{Lai2017}, incorporating additional operational constraints. Similarly, Andres et al. \cite{Andres2015} proposes scheduling maintenance based on flexible intervals rather than fixed mileage thresholds, directing train units to maintenance stations within allowable maintenance windows. As observed, the majority of these studies focus on routing railcars to various maintenance depots along their service paths. However, they typically do not address SLA requirements, maintenance depot capacity constraints, or breaches of maintenance intervals—critical elements that are explicitly considered in this study.


When examining the literature on stochastic failures, existing studies can be broadly categorized into three main domains: (i) track-only maintenance, (ii) joint track and train operations maintenance, and (iii) railcar maintenance specific to tram systems.


In the domain of track-only maintenance, one of the earliest contributions is presented by Consilvio et al. \cite{consilvio2016stochastic}, who incorporate probabilistic asset degradation to define maintenance deadlines and allocate teams under uncertain conditions. A similar approach is adopted by Malak et al. \cite{malak2023predictive}, who introduce chance constraints to account for degradation uncertainty. Although both studies employ stochastic modeling techniques, they do not consider operational factors such as train movements or maintenance depot capacity, which are addressed in the present work.

The second category, where maintenance activities affect both track and train operations, explores how uncertainty influences service continuity. For example, D’Ariano et al. \cite{d2019integrated} analyze the effect of stochastic maintenance durations on train timetables using a multi-objective optimization approach that minimizes deviations from planned schedules while clustering successive maintenance tasks. Bababeik et al. \cite{bababeik2022simultaneous} adopt chance-constrained modeling to capture uncertainty in track blockages and reduce delays, while Ji et al. \cite{ji2024optimization} employ two-stage stochastic programming to determine train departure times under uncertain maintenance durations. Unlike track-only models, these studies incorporate operational elements, such as maintenance capacity.

The third category focuses on railcar-specific maintenance, particularly within tram operations. An early study in this area is provided by Elhüseyni and Ünal \cite{elhuseyni2021integrated}, who model preventive, corrective, and out-of-service maintenance within a single-track environment while accounting for SLA violations and track capacity constraints. However, corrective maintenance decisions in their work are reactive and based on a rolling horizon framework. A different approach is taken by Del Castillo et al. \cite{del2023dynamic}, who develop a joint optimization model for train assignments and maintenance scheduling using RUL predictions derived from anomaly detection, aiming to balance railcar utilization and depot feasibility. Folco et al. \cite{folco2024rolling} extend this line of research by incorporating multiple maintenance cycles, modeling preventive and corrective costs based on failure probabilities and timing deviations, respectively, and applying a rolling horizon scheme. Finally, Liao et al. \cite{liao5026899scenario} address high-level maintenance planning under uncertainty in both failure timing and maintenance duration. Their two-stage stochastic programming approach determines start times and buffer windows across multiple scenarios, distinguishing itself by incorporating uncertainty into both strategic and operational maintenance layers.


The predictive maintenance approach based on failure scenarios is not limited to railcar maintenance, as similar examples can be found in power systems. One example is provided by Basciftci et al. \cite{Basciftci2018}, who study condition-based maintenance and operation scheduling of generators under unexpected failures using a two-stage stochastic model. This line of research is further extended by Okumuşoğlu et al. \cite{okumusoglu2022joint}, who incorporate transmission lines through an enhanced stochastic framework.


Table~\ref{tab:litable} presents a chronological overview of key studies on railway and railcar maintenance, classified by problem components, application areas, modeling assumptions, and solution approaches. Most studies incorporate train operations and maintenance capacity constraints, except for \cite{consilvio2016stochastic} and \cite{malak2023predictive}, which focus exclusively on asset degradation. Although deviations from target times are frequently analyzed, SLA and track violations are rarely integrated—explicitly addressed only in \cite{elhuseyni2021integrated} and in the present study. Preventive maintenance is widely considered, whereas corrective maintenance is typically handled in a reactive, rolling-horizon fashion, as in \cite{Lai2015,elhuseyni2021integrated,folco2024rolling}. With respect to uncertainty, predictive maintenance has been modeled through explicit stochastic approaches, employing either chance constraints or stochastic programming. Importantly, none of the reviewed works integrate all critical components—preventive and corrective maintenance, SLA and track violations, and stochastic failure scenarios—while simultaneously comparing deterministic MILP and stochastic programming outcomes within a unified framework based on SAA.

\begin{table}[h]
\centering
\caption{Overview of studies on railcar and railway maintenance (MILP: Mixed-integer Linear Programming, SP: Stochastic Programming, CC: Chance Constraint, MO: Multi-objective, SAA: Sample Average Approximation).}
\label{tab:litable}
\resizebox{\textwidth}{!}{%
\begin{tabular}{|c|c|c|c|c|c|c|c|c|c|c|c|}
\hline
\textbf{Papers}    & \multicolumn{8}{c|}{\textbf{Components}}                                                                                                                                        & \textbf{Application} & \textbf{Failure} & \textbf{Solution}   \\
\cline{2-9}
\textbf{ }          & {Train} & {Maintenance} & {SLA} & {Track} & {Deviation} & {Preventive} & {Corrective} & {Maintenance} & \textbf{Area}            & \textbf{Uncertainty}                   & \textbf{Method}                  \\
\textbf{}          & {operations} & {capacity} & {violation} & {violation} & {from target} & {maintenance} & {maintenance} & {action} & \textbf{}            & \textbf{}                   & \textbf{}                  \\
\hline
\cite{Lai2015} & yes                     & yes                     & no                & no                  & no                             & yes                       & yes                       & Preventive             & Railcar              & no                          & MILP                       \\
\cite{consilvio2016stochastic}      & no                      & no                      & no                & no                  & yes                            & yes                       & no                        & Predictive           & Railway              & yes                         & SP     \\
\cite{zhong2019rolling}          & yes                     & no                      & no                & no                  & no                             & yes                       & no                        & Preventive           & Railcar              & no                          & MILP                       \\
\cite{d2019integrated}       & yes                     & yes                     & no                & no                  & yes                            & yes                       & no                        & Preventive           & Railway              & yes                         & MO+MILP        \\
\cite{bababeik2022simultaneous}       & yes                     & yes                     & no                & no                  & no                             & yes                       & no                        & Predictive           & Railway              & yes                         & MILP+CC     \\
\cite{elhuseyni2021integrated} & yes                 & yes                     & yes               & yes                 & yes                            & yes                       & yes                       & Preventive             & Railcar              & no                          & MILP                       \\
\cite{malak2023predictive}      & no                      & no                      & no                & no                  & yes                            & yes                       & yes                       & Predictive           & Railway              & yes                         & MILP+CC     \\
\cite{del2023dynamic}  & yes                     & yes                     & no                & no                  & yes                            & yes                       & no                        & Preventive          & Railcar              & no                          & MILP                       \\
\cite{folco2024rolling}         & yes                     & yes                     & no                & no                  & yes                            & yes                       & yes                       & Preventive          & Railcar              & no                          & MILP                       \\
\cite{liao5026899scenario}         & yes                     & yes                     & no                & no                  & no                             & yes                       & no                        & Predictive           & Railcar              & yes                         & SP     \\
\cite{ji2024optimization}           & yes                     & yes                     & no                & no                  & no                             & yes                       & no                        & Predictive           & Railway              & yes                         & SP    \\
\textit{This work} & yes                     & yes                     & yes               & yes                 & no                            & yes                       & yes                       & Predictive           & Railcar              & yes                         & SP + SAA \\
\hline
\end{tabular}%
}
\end{table}

\subsection{Our Approach and Contributions}

In this paper, we develop a stochastic MILP with scenarios to solve the joint optimization of operations and maintenance scheduling for a railcar system operator that manages a fleet. Assuming that the railcar deterioration is modeled using the Weibull distribution (see recent studies~\citep{del2023dynamic, folco2024rolling} with a similar assumption), we estimate its RUL based on this distribution given the each railcar’s age (i.e., the time elapsed since last maintenance) to support a predictive maintenance approach. We generate failure scenarios based on the RUL distributions of each railcar in the fleet. Since the scenario space is exponentially large, we utilize an SAA approach to efficiently estimate a point estimate and obtain confidence intervals  for the optimal value. 
Finally, we compare the effectiveness of our proposed predictive scheduling framework against a preventive maintenance strategy inspired from~\cite{elhuseyni2021integrated}. Our computational experiments suggest that the proposed approach consistently yields lower total cost, SLA violations, track violations and corrective maintenance actions.


Our paper makes the following key contributions to the literature: 


\begin{itemize}
    \item Our work  significantly extends  \cite{elhuseyni2021integrated} by enabling predictive scheduling of maintenance activities, 
    incorporating two distinct SLA schedules, operational costs, variations in fleet size, and weekend-specific operational timetables.  
        \item To the best of our knowledge, we propose the first stochastic programming model developed in the literature handling preventive and corrective maintenance, SLA and track capacity violations simultaneously.
    \item Our paper is also the first study to apply the SAA method within a unified framework that includes both in-sample and out-of-sample evaluation in the context of railcar maintenance. 
        \item Through extensive computational experiments, we demonstrate the effectiveness of predictive maintenance compared to preventive maintenance.
\end{itemize}

The rest of the paper is structured as follows: In Section~\ref{s:problemFormulation}, we introduce the problem formulation, and present  deterministic optimization and stochastic programming approaches. In Section~\ref{s:solutionMethodology}, we provide the details of the key methodological components of our approach, namely,  scenario generation and  SAA procedures. In Section~\ref{sec:comp}, we present the results of our extensive computational experiments. Finally, in Section~\ref{s:conc}, we conclude our paper with final remarks.

\section{Problem Formulation} \label{s:problemFormulation}

In this section, we first describe our problem setting in Section~\ref{s:problemSetting}. Then, we present two different approaches based on a deterministic optimization model and a stochastic programming model in Sections~\ref{s:deterministicModel} and~\ref{s:stochasticModel}, respectively.

\subsection{Problem Setting} \label{s:problemSetting}

We consider a transportation  system with a fleet $\mathcal{J}$ of railcars  and a planning horizon with a set of time periods $\mathcal{T}$. 
We assume that the system operator is responsible for operational and maintenance decisions. 
Regarding the operational aspect, the primary goal of the system operator is to assign the railcars into the time periods so that passenger demand, given in the form of an SLA, is satisfied. In particular,  at least $SLA_t$ many railcars should be operational for each time period $t\in\mathcal{T}$, otherwise, a penalty of $C_s$  is applied for each unit of unsatisfied demand. In addition, each railcar incurs a cost of $C_o$ for each time period in which it is operational.

Since railcars in transportation systems are highly utilized, the system operator also needs to consider the maintenance aspect. We consider two types of maintenance actions: i) corrective maintenance, ii) predictive maintenance. If a railcar fails, an intermediate corrective maintenance action is taken, which takes $Y_c$ time periods and costs $C_c$ units. The system operator can proactively schedule a preventive maintenance as well, which takes $Y_p$ time periods ($Y_p<Y_c$) and costs $C_p$ units ($C_p<C_c$). We assume that there is a track capacity of $L$ allocated for maintenance actions and if this capacity is exceeded,  a penalty of $C_a$  is incurred for each additional maintenance activity.
%


A similar transportation system is also examined in~\cite{elhuseyni2021integrated}, however, the present study introduces several key distinctions. First, the Last-In-First-Out (LIFO) track parking constraints have been omitted to simplify the model. Second, unlike the referenced study, this work incorporates cost structures for each maintenance type, operational activity, and the utilization of additional tracks. These modifications aim to provide a more flexible and cost-aware maintenance planning framework.

\subsection{Deterministic Optimization Approach} \label{s:deterministicModel}

A commonly used approach in the maintenance planning literature is to schedule preventive maintenance actions according to manufacturer specifications. {Normally, a maintenance interval consists of begin and end points which are calculated according to a target age and a preferred tolerance, say $\pm10\%$ of the target age (see, e.g.,~\cite{elhuseyni2021integrated}). In this study, we assume that maintenance begin time of each railcar already overlaps with the planning period, so we discard it from a preferred maintenance interval.} Instead, we suppose that a preferred maintenance interval for railcar $j$ is given as $[D_j, \bar I_j]$  for $j \in\mathcal{J}$, where $D_j$  and $\bar I_j$ respectively represent the {due time} and the end time of the maintenance interval (we explain how we decide on these parameters in Section~\ref{sec:data}).  We assume that a unit penalty is applied for each unit of tardiness and earliness with respect to  due time $D_j$, and  $C_v>1$ units of penalty are incurred for each unit of maintenance violation beyond $\bar I_j$. The rationale behind this parameter selection is that after the end time of the preferred maintenance interval $\bar I_j$, the deterioration of the reliability increases dramatically, hence, this occurrence is heavily penalized~\cite{elhuseyni2021integrated}.

Having defined problem components, we now define our preventive maintenance problem. 
We tabulate the decision variables of our deterministic optimization model in Table~\ref{tab:decVarDeter}.
\begin{table}[H]
\centering
\small
\caption{Decision variables of the deterministic optimization model.}
\label{tab:decVarDeter}
\begin{tabular}{|cl}
\hline
\multicolumn{2}{|l|}{$\textbf{Decision Variables}$}     \\ \hline
    $z_{jt}$ & \multicolumn{1}{l|}{1 if railcar $j$ enters maintenance in period $t$, and 0 otherwise.}\\
    $p_{j}$ & \multicolumn{1}{l|}{1 if a maintenance is not scheduled for railcar $j$ within the planning horizon, and 0 otherwise.} \\
    $\eta_{jt}$ & \multicolumn{1}{l|}{1 if railcar $j$ is operational in period $t$, and 0 otherwise.} \\  
    $\sigma_t$ &  \multicolumn{1}{l|}{Number of railcars below $SLA_t$.}\\
    $\gamma_t$ & \multicolumn{1}{l|}{Additional maintenance capacity added in period $t$.}\\
    $v_j$ &  \multicolumn{1}{l|}{Number of time units the maintenance interval of job $j$ is violated.}\\
    $s_{j}$ & \multicolumn{1}{l|}{Start time of job $j$.} \\
    $e_{j}$ & \multicolumn{1}{l|}{Earliness  of job $j$.} \\
    $t_j$ & \multicolumn{1}{l|}{Tardiness of job $j$.} \\

\hline
\end{tabular}
\end{table}

We now present the deterministic optimization  model adopted from~\cite{elhuseyni2021integrated} as follows:
\begin{subequations}   \label{eq:deterministicModelFormulation}
	\begin{align}
		\min &\hspace{.5em} 
        C_o \sum_{j \in \mathcal{J}} \sum_{t \in \mathcal{T}} \eta_{jt}  + C_s \sum_{t \in \mathcal{T}} \sigma_t   + C_{a} \sum_{t \in \mathcal{T}}  \gamma_{t}     +  C_p \sum_{j \in \mathcal{J}} \sum_{t \in \mathcal{T}}  z_{jt} + \sum_{j \in \mathcal{J}} \left( t_j+e_j + C_{v} v_j\right) 
         \label{objfuncDeter}\\
  \mathrm{s.t.} 
  		&\hspace{0.5em} \sum_{t \in \mathcal{T}} z_{jt} + p_j = 1  & j \in \mathcal{J}	\label{jobassignDeter}\\
  		&\hspace{0.5em} t_{j} \ge s_j + \vert \mathcal{T} \vert - D_j - \vert \mathcal{T} \vert(1-p_j) & j \in \mathcal{J} \label{tardiness}\\ 
		&\hspace{0.5em} e_{j} \ge D_j-s_j - \vert \mathcal{T} \vert p_j& j \in \mathcal{J} \label{earlinesseither} \\
		&\hspace{0.5em} t_{j} - v_j \le \bar I_j-D_j& j \in \mathcal{J}  \label{tardUB}\\ 
        &\hspace{0.5em} \sum_{t \in \mathcal{T}} t z_{jt}  = s_j  & j  \in \mathcal{J}	\label{jobbegin}\\
        &\hspace{0.5em} \sum_{j \in J} \eta_{jt} + \sigma_t\ge SLA_t  & t \in \mathcal{T} \label{SLAsatisfDeter} \\
    &\hspace{0.5em} \eta_{jt} \le 1-\sum_{e=0}^{Y_p - 1} z_{j, t-e} & j \in \mathcal{J}, t \in \mathcal{T} \label{availib}  \\ 
    &\hspace{0.5em} \sum_{j \in \mathcal{J}} \sum_{e=0}^{Y_p - 1} z_{j, t-e} \le L + \gamma_{t}  & t \in \mathcal{T} \label{capacityDeter} \\  
%
%
            &\hspace{0.5em} t_j, e_j, s_j, v_j \ge 0   & j \in \mathcal{J} \label{tardvar} \\
            &\hspace{0.5em} \sigma_t, \gamma_{t} \ge 0 & t \in \mathcal{T} \\
            &\hspace{0.5em} z_{jt}, p_j, \eta_{jt} \in \{0,1\}  & j \in \mathcal{J}, t \in \mathcal{T} . \label{binaryvarsDeter}
		\end{align}
\end{subequations}
Here, the objective function~\eqref{objfuncDeter} minimizes the total cost with five components:  i) the operational cost,  ii) the SLA violation penalty, iii) the additional maintenance cost exceeding the track capacity, iv) the preventive maintenance cost, and v) the deviations from a ``preferred" maintenance interval. 
Constraint~\eqref{jobassignDeter} guarantees that each job is assigned a predictive maintenance period or is postponed. Constraints~\eqref{tardiness} and \eqref{earlinesseither} quantify the tardiness and the earliness, respectively, whereas constraint~\eqref{tardUB} computes the maintenance interval violations. Constraint~\eqref{jobbegin} assigns the start time of jobs while constraint~\eqref{SLAsatisfDeter} computes the SLA violations. Constraint~\eqref{availib} makes sure that a railcar is not operational when it is under preventive maintenance. Constraint~\eqref{capacityDeter} is the track capacity constraint, which also keeps track of additional maintenance activities. 
Finally, constraints~\eqref{tardvar}-\eqref{binaryvarsDeter} are variable domain restrictions. 

\subsection{Stochastic Programming Approach} \label{s:stochasticModel}

Although deterministic optimization approaches are widely used in maintenance planning, they do not explicitly consider the unexpected failures and the resulting {out-of-service states}. For example, the MILP model~\eqref{eq:deterministicModelFormulation} encourages the decision-maker to schedule preventive maintenance until the due time $D_j$, which can be set according to producer specifications. However, even under such a policy, unexpected failures can occur, corrective maintenance action needs to be taken and system disruptions can happen. In this section, we propose a stochastic programming model, which explicitly considers unexpected failures through scenarios in order to manage such disruptions in a cost-effective manner.


Let us denote the failure time of railcar $j$ with the random variable $\Xi_j$. We assume that we are given a finite set of scenarios denoted by $\mathcal{K}$, and denote the probability of scenario $k\in\mathcal{K}$ as $\Pi^k$. We will denote the {failure time of railcar $j$ in scenario $k$} as $\Xi^k_j$ (we explain how we generate our scenarios in Section~\ref{sec:data}). 
We tabulate the decision variables of our two-stage stochastic programming model in Table~\ref{tab:decVarStoch}.

\begin{table}[H]
\centering
\small
\caption{Decision variables of the stochastic programming model.}
\label{tab:decVarStoch}
\begin{tabular}{|cl}
\hline
\multicolumn{2}{|l|}{$\textbf{Decision Variables}$}     \\ \hline
    $z_{jt}$ & \multicolumn{1}{l|}{1 if railcar $j$ enters maintenance in period $t$, and 0 otherwise.}\\
    $p_{j}$ & \multicolumn{1}{l|}{1 if a maintenance is not scheduled for railcar $j$ within the planning horizon, and 0 otherwise.} \\
    $\eta^k_{jt}$ & \multicolumn{1}{l|}{1 if railcar $j$ is operational in period $t$ in scenario $k$, and 0 otherwise.} \\ 
        {$m^k_{jt}$} & \multicolumn{1}{l|}{{1 if railcar $j$ is under maintenance in period $t$ in scenario $k$, and 0 otherwise.}}   \\
    $\sigma_t^{k}$ &  \multicolumn{1}{l|}{Number of railcars below $SLA_t$ in scenario $k$.}\\
    $\gamma^k_{t}$ & \multicolumn{1}{l|}{Additional maintenance capacity added in period $t$ in scenario $k$.}    \\
\hline
\end{tabular}
\end{table}
We note that maintenance variables $z_{jt}$ and $p_j$ are \textit{first-stage} decisions, meaning that they are decided before any uncertainty is realized. On the other hand, the other variables correspond to \textit{second-stage} decisions, which are a function of the uncertainty observed (hence, the $k$ superscript).
Now, we are ready to present the stochastic programming model formulation, which is partially inspired from maintenance scheduling studies ~\cite{Basciftci2018,okumusoglu2022joint} from the power systems literature, as follows:

\begin{subequations}   \label{eq:stochasticModelFormulation}
	\begin{align}
		\min &\hspace{0.5em}   C_o \sum_{k \in \mathcal{K}}\Pi^k  \sum_{j \in \mathcal{J}} \sum_{t \in \mathcal{T}} \eta^k_{jt} +  C_s \sum_{k \in \mathcal{K}}\Pi^k \sum_{t \in \mathcal{T}} \sigma_t^{k}   \nonumber \\
          \hspace{1.5em} & +   C_{a} \sum_{k \in \mathcal{K}}\Pi^k \sum_{t \in \mathcal{T}}  \gamma^k_{t}  +  
          C_p  \sum_{k \in \mathcal{K}}\Pi^k\sum_{j \in \mathcal{J}} \sum_{t=1}^{\Xi^k_j-1} z_{jt}  \nonumber \\
          \hspace{1.5em}   &+  C_c \sum_{k \in \mathcal{K}}\Pi^k \sum_{j \in \mathcal{J}} \sum_{t=\Xi^k_j}^{\vert \mathcal{T} \vert} (z_{jt}+p_j) 
        \label{objfuncStoch}\\
		\mathrm{s.t.}   	
    %
&\hspace{0.5em} \sum_{t \in \mathcal{T}} z_{jt} + p_{j} =1  & j \in \mathcal{J}   \label{jobassignStoch}\\
%
    &\hspace{0.5em} \sum_{j \in J} \eta^k_{jt} + \sigma_t^{k}\ge SLA_t  & t \in \mathcal{T},  k \in \mathcal{K} \label{SLAsatisfStoch} \\  
    &\hspace{0.5em} m_{jt}^k = \sum_{e=0}^{\min\{Y_p,t\}-1} z_{j,t-e} & j \in \mathcal{J}, k \in \mathcal{K} , t=1,\dots,\Xi_j^k-1 \label{prevavailib}  \\ 
    & \hspace{0.5em}   m_{jt}^k = 1-\sum_{t'=1}^{\min\{t-Y_p, \Xi_j^k-1 \}} z_{j,t'} \quad &j\in\mathcal{J}, k\in\mathcal{K} , 
    t=\Xi_j^k, \dots, \min\{\Xi_j^k+Y_c-1, |\mathcal{T}|\} 
    \label{coravailib}\\
    & \hspace{0.5em}  m_{jt}^k = 0 \quad &j\in\mathcal{J}, k\in\mathcal{K} , t= \Xi_j^k+Y_c, \dots, |\mathcal{T}| \label{nonavailib}\\
        & \hspace{0.5em}  \eta_{jt}^k \le 1- m_{jt}^k   \quad & j\in\mathcal{J}, k\in\mathcal{K} , t\in\mathcal{T} \label{opermaintlink}  \\
        & \hspace{0.5em}  \sum_{j \in J} m_{jt}^k \le L + \gamma_t^k  \quad &  k\in\mathcal{K} , t\in\mathcal{T} 
        \label{capacityStoch} \\
    %
    &\hspace{0.5em} \sigma_t^{k}, \gamma^k_{t} \ge 0 & t \in \mathcal{T}, k \in \mathcal{K} \label{contVarStoch} \\
            &\hspace{0.5em}  z_{jt}, p_j, \eta^k_{jt},  m^k_{jt} \in \{0,1\}  & j \in \mathcal{J}, t \in \mathcal{T}, k \in \mathcal{K}\label{binVarStoch}.   
    %
		\end{align}
\end{subequations}
Here, the objective function~\eqref{objfuncStoch} minimizes the \textit{expected} total cost with five components: i) the expected operational cost, ii) the expected SLA violation penalty, iii) the expected additional maintenance cost exceeding the track capacity, iv) the expected preventive maintenance cost, and v) the expected corrective maintenance cost. 
Constraint~\eqref{jobassignStoch} guarantees that each job is assigned at most one maintenance action. 
Constraint~\eqref{SLAsatisfStoch} computes the SLA violations for each scenario. 
Constraints~\eqref{prevavailib}-\eqref{nonavailib} make sure that a railcar is in maintenance state based on its scheduled maintenance time and the failure time for each scenario. 
{To be more specific, if a maintenance decision is given in the first $Y_p$ periods, it becomes preventive with regard to constraint~\eqref{prevavailib}. Starting period $\Xi_j^k$, constraint~\eqref{coravailib} stipulates that either there is a continuing PM job, or the job becomes corrective due to the lack of maintenance decision before period $\Xi_j^k$. A job cannot be at maintenance state after $\Xi_j^k+Y_c$ periods with regard to \eqref{nonavailib}.} 
{Constraint ~\eqref{opermaintlink} links maintenance state to operational state.} Constraint~\eqref{capacityStoch} is the track capacity constraint, which also keeps track of additional maintenance activities for each scenario. Finally, constraints~\eqref{contVarStoch}-\eqref{binVarStoch} are variable domain restrictions.

\section{Solution Methodology} \label{s:solutionMethodology}

In this section, we present the main components of our solution methodology. In particular, we first describe the data generation procedure in Section~\ref{sec:data}. Then, we introduce SAA approach in Sections~\ref{sec:inSampleOpt} and~\ref{sec:outOfSampleEval} that respectively focus on the in-sample optimization step and the out-of-sample evaluation step.



\subsection{Data/Scenario Generation}
\label{sec:data}

 
An important component in our analysis is the lifetime distribution of the railcars. In our paper, we assume that the lifetime of a railcar follows a Weibull distribution with a scale parameter $\alpha$ and a shape parameter $\beta$ that has the following probability density function (pdf):
\[
 f(y) := \frac{\beta}{\alpha} \left(\frac{y}{\alpha}\right)^{\beta-1} e^{- (y/\alpha)^{\beta}  }, \ y\ge0.
\]
We assume that the \textit{initial age} of a railcar, defined as the mileage after the last maintenance action at the beginning of planning, is given as $y'$. Then, we can compute the RUL distribution as follows~\cite{nelson1982applied}:
\[
f(y|y') = \frac{\beta}{\alpha} \left(\frac{y}{\alpha}\right)^{\beta-1} e^{(y'/\alpha)^{\beta} - (y/\alpha)^{\beta} }, \ y\ge y'.
\]
In particular, the conditional probability of a railcar not failing before time $y$ given an initial age $y'$ can be computed as
\begin{equation}\label{eq:WeibullConditionalCDF}
    F(y|y') = 1-e^{  (y'/\alpha)^{\beta}-(y/\alpha)^{\beta} }, \ y \ge y'.
\end{equation}
We present an illustrative example in Figure~\ref{fig:fleetsTogether} for three different initial age values. Here, the reference point $x=0$ is the beginning of the planning at which the initial age $y'$ is given. 

\newcommand\youngAge{10}
\newcommand\mediumAge{30}
\newcommand\oldAge{50}
\newcommand\WeibullAlpha{50}
\newcommand\WeibullBeta{5}

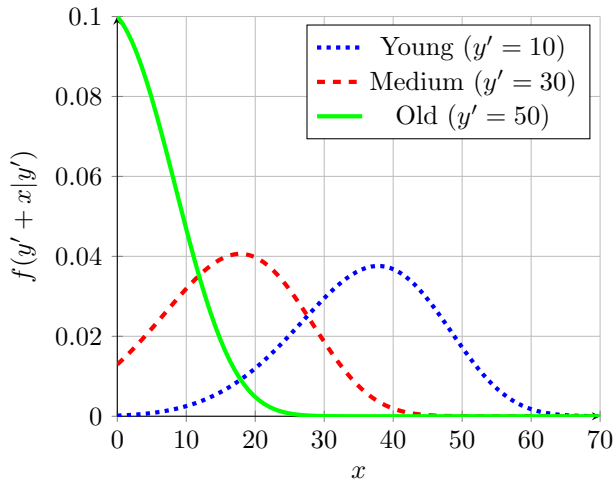
\begin{figure}[htbp]
    \centering
        \begin{tikzpicture}
        \begin{axis}[
            axis lines=left,
            xlabel={$x$},
            ylabel={$f(y' + x |  y')$},
            domain=0:70,
            ymax=0.1,
            samples=100,
            width=8cm,
            grid=major,
            yticklabel style={/pgf/number format/.cd,fixed,precision=3},
        ]
            \addplot[blue, ultra thick, dotted] {(\WeibullBeta/\WeibullAlpha)*(((x+\youngAge)/\WeibullAlpha)^(\WeibullBeta-1)*exp((\youngAge/\WeibullAlpha)^\WeibullBeta - ((x+\youngAge)/\WeibullAlpha)^\WeibullBeta)};

            \addlegendentry{Young ($y'=\youngAge)$}; 
            
            \addplot[ red, ultra thick, dashed] {(\WeibullBeta/\WeibullAlpha)*(((x+\mediumAge)/\WeibullAlpha)^(\WeibullBeta-1)*exp((\mediumAge/\WeibullAlpha)^\WeibullBeta - ((x+\mediumAge)/\WeibullAlpha)^\WeibullBeta)};
                        
            \addlegendentry{Medium ($y'=\mediumAge)$};

            \addplot[green, ultra thick] {(\WeibullBeta/\WeibullAlpha)*(((x+\oldAge)/\WeibullAlpha)^(\WeibullBeta-1)*exp((\oldAge/\WeibullAlpha)^\WeibullBeta - ((x+\oldAge)/\WeibullAlpha)^\WeibullBeta)};

            \addlegendentry{Old ($y'=\oldAge)$};

        \end{axis}
        \end{tikzpicture}

    \caption{Weibull remaining lifetime distributions (pdf) for different initial age values ($\alpha=\WeibullAlpha$, $\beta=\WeibullBeta$).}
    \label{fig:fleetsTogether}
\end{figure}
In our experiments, we generate failure scenarios by sampling from the RUL distribution given in~\eqref{eq:WeibullConditionalCDF} via the inverse transform method. In particular, given $U=\text{Unif}(0,1)$, we set the failure time as
\[ \Xi = \min\left\{ \left \lceil \alpha [ \left( {y'}/{\alpha} \right)^{\beta} - \ln(1 - U) ]^{\frac{1}{\beta}} \right\rceil - y' , \ |\mathcal{T}|+1 \right\}. \] 
Here, we set the failure time as $|\mathcal{T}|+1$ if the sampled time is outside the planning horizon. We repeat this procedure for each railcar $j$ and scenario $k$ to obtain the failure times $\Xi_j^k$. 
We use this scenario generation procedure while constructing the stochastic programming model~\eqref{eq:stochasticModelFormulation} as well as in the out-of-sample evaluation as described in Section~\ref{sec:outOfSampleEval}.

Two important inputs in the deterministic optimization model~\eqref{eq:deterministicModelFormulation} are the due date and the end time of the maintenance interval. We again use the Weibull distribution CDF, given as $F(y) = 1-e^{-(y/\alpha)^{\beta}}$, to decide these parameters. In particular, given a reliability level $R \in (0,1)$, we first compute the time time at which the railcar fails with probability $R$ as $\tau_R := \left\lceil \alpha \cdot \left( -\ln(1 - R) \right)^{1/\beta} \right\rceil$. Then, we compute a \textit{maintenance interval half-length}, denoted by, $\eta_R := \lceil 0.1\tau_R \rceil$. Finally, given a railcar with an initial age of $y'$, we set the due date as $\lceil \tau_R+\frac12 \eta_R \rceil$ and the end time of the maintenance interval as $\tau_R+ \eta_R $ 
according to \cite{elhuseyni2021integrated}.
For the example described in Figure~\ref{fig:fleetsTogether}, if we set $R=0.8$, the maintenance intervals, whose end points are the due date $D_j$ and the end date $\bar I_j$, are calculated as $[48,51]$, $[28,31]$ and $[8,11]$ for railcars with a young, medium and old initial age, respectively.
Notice that with larger values of $R$, the maintenance interval shifts to the right, which increases the utilization of the railcar and the risk of failure at the same time. For instance, if we set $R=0.9$, the maintenance intervals mentioned before become $[53,56]$, $[33,36]$ and $[13,16]$ for railcars with a young, medium and old initial age, respectively.



\subsection{In-Sample Optimization}
\label{sec:inSampleOpt}

In the deterministic optimization approach, we simply solve the MILP~\eqref{eq:deterministicModelFormulation}. 
Let us denote an optimal maintenance solution obtained this way as $z^D$.

In the stochastic programming approach, we adopt the SAA approach. In particular, we first solve the MILP~\eqref{eq:stochasticModelFormulation} with a finite sample $\mathcal{K}$ with $|\mathcal{K}|=N$. We then  repeat this procedure $M$ times due to this randomness, and record their optimal values as $\nu_m$ and optimal maintenance solutions. Then, we set $z^S$ as the maintenance solution that gives the smallest  objective function value $\nu_m$ value among the $M$ replications. Different from the deterministic optimization approach, we can obtain a confidence interval (CI) which likely contains a lower bound estimate for the unknown optimal value. For this purpose, we first  calculate the mean and variance estimates of the true lower bound estimate as follows:
\[\hat \mu_L = \frac{1}{M} \sum_{m=1}^{M} \hat \nu_m \quad \text{ and } \quad \hat \sigma^2_{L} = \frac{1}{M(M-1)} \sum_{m = 1}^{M} \big( \nu_m - \hat \mu_L\big)^2.\]
Afterwards, we construct the approximate $(1-\theta)$ level CI for the lower bound estimate as $ \hat \mu_L \pm t_{\theta/2, M-1} \hat \sigma_L$.

\subsection{Out-of-Sample Evaluation}
\label{sec:outOfSampleEval}

Given a maintenance decision~$z^*$ (which is~$z^D$ for the deterministic approach or~$z^S$ for the stochastic approach), we would like to evaluate (or test) its success in a fresh sample. For this purpose, we generate a new sample with $N' \gg N$ and compute $\omega_n$ for $n=1,\dots,N'$ by running Algorithm~\ref{alg:second-stage} with the fixed maintenance decision~$z^*$ for each sample~$n$ separately. 

\begin{algorithm}
\caption{The second-stage objection function value calculation.}
\label{alg:second-stage}
    \begin{algorithmic}
    \REQUIRE A first-stage decision $ z_{jt}^*$ and a scenario index $n$.
    \ENSURE The second-stage objection function value $\omega_n$.  
    \STATE Initialize $\omega_n = 0$.
    \FOR{$j \in J$}
        \STATE Compute the time the maintenance is scheduled for railcar $j$ as 
        $\hat m_j = \sum_{t\in\mathcal{T}} tz_{jt}^*$. 
        %
        \IF{$1 \le \hat m_j < \Xi_j^n$}
        \STATE Update $\omega_n = \omega_n + C_p$.
        \STATE Set the times the railcar is under preventive maintenance as $U_j^n =\{ \hat m_j, \hat m_j + 1 , \dots, \min\{ \hat m_j + Y_p-1, |\mathcal{T}| \}$.
        \ENDIF
        \IF{$\Xi_j^n \le \hat m_j \le |\mathcal{T}|$ or {($ \hat m_j = 0$ and $\Xi_j^n \le |\mathcal{T}|$)}}
        \STATE Update $\omega_n = \omega_n + C_c$.
        \STATE Set the times the railcar is under corrective maintenance as $U_j^n =\{ \Xi_j^n, \Xi_j^n + 1 , \dots, \min\{ \Xi_j^n + Y_c-1, |\mathcal{T}| \}$.
        \ENDIF
        \ENDFOR
    \FOR{$t \in T$}
        \STATE Compute the number of railcars under maintenance in period $t$ as $M_t^n = \sum_{j\in J} \boldsymbol{1}(t \in U_j^n)$.
        \STATE Compute
        the number of available railcars in period $t$ as $A_t^n = |\mathcal{J}|-M_t^n$.
        \STATE Update $\omega_n = \omega_n + C_o \min\{A_t^n, SLA_t\} + C_s \max\{SLA_t- A_t^n, 0\}  + C_a \max\{M_t^n - L , 0\}$.
    \ENDFOR
    
    \end{algorithmic}
\end{algorithm}

After computing these objective function values, we can obtain a  CI  which likely contains an upper bound estimate for the unknown optimal value. For this purpose, we first  calculate the mean and variance estimates of the true upper bound estimate as follows:
\[\hat \mu_U = \frac{1}{N'} \sum_{n=1}^{N'} \hat \omega_n \quad \text{ and } \quad \hat \sigma^2_{U} = \frac{1}{N'(N'-1)} \sum_{n = 1}^{N'} \big( \omega_n - \hat \mu_U\big)^2.\]
Afterwards, we construct the approximate $(1-\theta)$ level CI for the upper bound estimate as $  \hat \mu_U \pm z_{\theta/2} \hat \sigma_U$.

\section{Computational Experiments}
\label{sec:comp}


In this section, we discuss the results of our computational experiments. 
After providing the computational setting in Section~\ref{sec:compSetting}, we present the results of our out-of-sample analysis that compares the stochastic  and deterministic models in Section~\ref{sec:compOutOfSample}. Then, we carry out a confidence interval analysis to justify the reliability of the parameter selections in Section~\ref{sec:compCI}. Finally, we provide a detailed solution comparison in Section~\ref{sec:detailedSolComp} to illustrate the key differences in the optimal solutions of stochastic and deterministic models.

\subsection{Setting}
\label{sec:compSetting}

In reference to Section~\ref{s:problemFormulation}, we have chosen the problem parameters as follows: {$Y_p=3$, $Y_c=5$, $C_o=2$, $C_{v} = 2 C_o$, {$C_c = 20 C_o$}, $C_p=C_a=C_c/5$}. Since the relationship between the SLA violation penalty $C_s$ and the corrective cost $C_c$ is quite influential on the outcomes, we carry out a parametric analysis as follows:
\begin{itemize}
    \item Case 1: $C_s = \frac32 C_c $. In this case, the SLA violation penalty is higher than the corrective cost.
    \item Case 2: $C_s = C_c $. In this case, the SLA violation penalty and the corrective cost are equal.
    \item Case 3: $C_s = \frac23 C_c $. In this case, the SLA violation penalty is lower than the corrective cost.
\end{itemize}
Under each of these cases, we test the effect of two other parameters: i) SLA of metro lines, ii) age distribution of the fleet. We take the planning horizon as seven days where each day consists of four time periods: Night, Morning, Afternoon, and Evening. In our case study, we use the schedules of two metro lines from Istanbul~\cite{metroistanbul2025}: M1B and M4.  
{We tally the number of services at each period and Table~\ref{tab:slasched} presents the resulting schedule. The maximum number of railcars in M1B and M4 are 51 and 56, respectively. We assume that the system works with 10\% of spare railcars, resulting in $n=57$ and $n=74$ railcars as fleet size for M1B and M4, respectively. We set the track size 20\% of the corresponding fleet size, leading to $L=12$ and $L=15$ for M1B and M4, respectively.}

\begin{table}[htbp]
\small
\centering
\caption{SLA with respect to different lines and daytime.}
\label{tab:slasched}
\begin{tabular}{|c|ccc|ccc|} \hline
          &       & M1B      &        &       & M4       &        \\ \hline
Daytime   & Weekdays & Saturday & Sunday & Weekdays & Saturday & Sunday \\ \hline
Night  [00:00-06:00]     & 0     & 12       & 12     & 0     & 11       & 11 \\
Morning [06:00-12:00]   & 45    & 41       & 40     & 68    & 55       & 48     \\
Afternoon [12:00-18:00] & 50    & 51       & 43     & 66    & 66       & 66     \\
Evening [18:00-00:00]   & 49    & 47       & 41     & 59    & 56       & 54     \\ \hline
\end{tabular}
\end{table}

The lifetime  of railcars is assumed to come from a Weibull distribution with a scale parameter $\alpha=\WeibullAlpha$ and a shape parameter $\beta=\WeibullBeta$.  
As for the initial age distribution, we consider a `young' fleet, an `old' fleet or a `mixed' one.  A young fleet refers to a fleet in which the initial ages are sampled from {$\text{Unif}(10, 20)$} whereas an old fleet refers   to a fleet in which the initial ages are sampled from {$\text{Unif}(50, 60)$}. In a mixed fleet, the initial ages are sampled from {$\text{Unif}(10, 60)$}.

We compare the following solution methods: 
\begin{itemize}
    \item \texttt{Stochastic Model}: Stochastic Programming Formulation~\eqref{eq:stochasticModelFormulation} with sample sizes {$|\mathcal{K}|=N=150$ and $M=5$.} 
    \item \texttt{Strict Deterministic Model}: Deterministic Optimization Formulation~\eqref{eq:deterministicModelFormulation}, where the maintenance interval is constructed with $R=0.8$ as described in Section~\ref{sec:data}.
    \item \texttt{Relaxed Deterministic Model}: Deterministic Optimization Formulation~\eqref{eq:deterministicModelFormulation}, where the maintenance interval is constructed with $R=0.9$.
\end{itemize}


We choose   $N'=1000$ in the out-of-sample evaluation. 
%
We have five key performance indicators (KPI), which are averaged over out-of-sample evaluations, for each case as listed below: 
\begin{itemize}
    \item \#Prev: The   number of ongoing preventive jobs, averaged over the periods.
    \item \#Cor: The average number of ongoing corrective jobs, averaged over the periods.
    \item SLA-v: The  {average} SLA violation, averaged over the periods.
    \item Track-v: The {average} track violation, averaged over the periods.
    \item Cost: The total cost (including maintenance and operational costs) estimate $\hat \mu_U$ obtained from the out-of-sample evaluation described in Subsection~\ref{sec:outOfSampleEval}.  
\end{itemize}
We also report the total CPU time of in-sample optimization  in seconds for each experiment.
We conduct our computational experiments using C++ calling CPLEX 22.1.1 solver running on 32 threads on a PC with Intel(R) Xeon(R) Silver 4210R CPU having two processors and 64 GB RAM. A time limit of three hours 
is enforced for each MILP.

\subsection{Out-of-Sample Analysis}
\label{sec:compOutOfSample}

We report the results of our computational study in Tables~\ref{tab:M1B} and \ref{tab:M4} for Lines M1B and M4, respectively.  {Note that the testing time is negligible, hence, it is not shown in these tables.}

\begin{landscape}

\begin{table}[h]
\centering
\caption{Computational results for Line M1B.}
\label{tab:M1B}
\footnotesize
\begin{tabular}{|c|c|cccccc|cccccc|cccccc|}
\hline
\multicolumn{2}{|c|}{\textbf{Setting}} 
& \multicolumn{6}{c|}{\texttt{Stochastic Model}} 
& \multicolumn{6}{c|}{\texttt{Strict Deterministic}} 
& \multicolumn{6}{c|}{\texttt{Relaxed Deterministic}} \\
\hline
Fleet & Case 
& \#Prev & \#Cor & SLA\_v & Track\_v & Cost & Time(s) 
& \#Prev & \#Cor & SLA\_v & Track\_v  & Cost & Time(s) 
& \#Prev & \#Cor & SLA\_v & Track\_v  & Cost& Time(s)  \\
\hline
     \multirow{3}{*}{Young} &      1 &       5.17 &       0.52 &       0.02 &       0.01 &    2559.04 &      12038 &       0.00 &       3.19 &       0.05 &       0.00 &    2950.52 &          0 &       0.00 &       3.19 &       0.05 &       0.00 &    2950.52 &          0 \\

       &      2 &       5.17 &       0.51 &       0.02 &       0.01 &    2549.76 &      10100 &       0.00 &       3.19 &       0.05 &       0.00 &    2925.28 &          0 &       0.00 &       3.19 &       0.05 &       0.00 &    2925.28 &          0 \\

       &      3 &       5.30 &       0.47 &       0.03 &       0.01 &    2544.51 &       1350 &       0.00 &       3.19 &       0.05 &       0.00 &    2907.61 &          0 &       0.00 &       3.19 &       0.05 &       0.00 &    2907.61 &          0 \\
\hline
     \multirow{3}{*}{Mixed} &      1 &       4.83 &       1.72 &       0.21 &       0.00 &    3117.47 &      43947 &       0.83 &       6.01 &       0.40 &       0.00 &    4154.50 &          0 &       0.24 &       6.84 &       0.53 &       0.00 &    4531.08 &          0 \\

       &      2 &       4.93 &       1.56 &       0.25 &       0.00 &    3009.48 &      41878 &       0.83 &       6.01 &       0.40 &       0.00 &    3931.90 &          0 &       0.24 &       6.84 &       0.53 &       0.00 &    4234.96 &          0 \\

       &      3 &       5.06 &       1.33 &       0.28 &       0.01 &    2897.40 &       7927 &       0.83 &       6.01 &       0.40 &       0.00 &    3776.08 &          0 &       0.24 &       6.84 &       0.53 &       0.00 &    4027.67 &          0 \\
\hline
       \multirow{3}{*}{Old} &      1 &       4.65 &       2.42 &       2.39 &       0.00 &    7429.13 &       9226 &       2.61 &       5.83 &       3.02 &       0.00 &    9077.54 &          1 &       1.97 &       6.90 &       3.25 &       0.00 &    9713.70 &          2 \\

         &      2 &       4.86 &       2.08 &       2.41 &       0.00 &    6080.62 &      10181 &       2.62 &       5.81 &       3.02 &       0.00 &    7384.48 &          4 &       1.95 &       6.93 &       3.26 &       0.00 &    7906.12 &          3 \\

         &      3 &       4.87 &       2.06 &       2.42 &       0.00 &    5132.71 &       3158 &       2.62 &       5.81 &       3.02 &       0.00 &    6200.70 &          1 &       1.95 &       6.93 &       3.26 &       0.00 &    6636.02 &          9 \\

\hline

\end{tabular}
\end{table}

\begin{table}[h]
\centering
\caption{Computational results for Line M4.}
\label{tab:M4}
\footnotesize
\begin{tabular}{|c|c|cccccc|cccccc|cccccc|}
\hline
\multicolumn{2}{|c|}{\textbf{Setting}} 
& \multicolumn{6}{c|}{\texttt{Stochastic Model}} 
& \multicolumn{6}{c|}{\texttt{Strict Deterministic}} 
& \multicolumn{6}{c|}{\texttt{Relaxed Deterministic}} \\
\hline
Fleet & Case 
& \#Prev & \#Cor & SLA\_v & Track\_v & Cost & Time(s) 
& \#Prev & \#Cor & SLA\_v & Track\_v  & Cost & Time(s) 
& \#Prev & \#Cor & SLA\_v & Track\_v  & Cost& Time(s)  \\
\hline

     \multirow{3}{*}{Young}  &      1 &       5.23 &       1.16 &       0.03 &       0.01 &    3400.82 &       4711 &       0.00 &       3.64 &       0.10 &       0.00 &    3816.50 &          0 &       0.00 &       3.64 &       0.10 &       0.00 &    3816.50 &          0 \\

       &      2 &       5.47 &       1.04 &       0.03 &       0.01 &    3381.70 &       4226 &       0.00 &       3.64 &       0.10 &       0.00 &    3759.86 &          0 &       0.00 &       3.64 &       0.10 &       0.00 &    3759.86 &          0 \\

       &      3 &       5.65 &       1.02 &       0.04 &       0.01 &    3376.36 &       1789 &       0.00 &       3.64 &       0.10 &       0.00 &    3720.21 &          0 &       0.00 &       3.64 &       0.10 &       0.00 &    3720.21 &          0 \\
\hline
      \multirow{3}{*}{Mixed} &      1 &       5.15 &       3.23 &       1.08 &       0.01 &    5593.81 &      54030 &       1.65 &       7.07 &       1.25 &       0.00 &    6523.81 &          0 &       0.57 &       8.78 &       1.57 &       0.00 &    7414.81 &          0 \\

       &      2 &       5.54 &       2.84 &       1.12 &       0.01 &    4970.76 &      40736 &       1.67 &       7.06 &       1.25 &       0.00 &    5820.86 &          0 &       0.57 &       8.78 &       1.57 &       0.00 &    6533.17 &          0 \\

       &      3 &       5.94 &       2.32 &       1.22 &       0.01 &    4512.51 &      12557 &       1.65 &       7.07 &       1.25 &       0.00 &    5329.80 &          0 &       0.57 &       8.78 &       1.57 &       0.00 &    5916.02 &          0 \\
\hline
        \multirow{3}{*}{Old} &      1 &       6.07 &       3.10 &       3.73 &       0.01 &   10737.90 &      26353 &       2.45 &       9.14 &       4.99 &       0.00 &   14009.60 &          1 &       1.51 &      10.70 &       5.30 &       0.00 &   14913.80 &          9 \\

         &      2 &       6.27 &       2.76 &       3.76 &       0.01 &    8644.96 &       9899 &       2.77 &       8.60 &       4.87 &       0.00 &   10931.70 &          1 &       1.92 &      10.02 &       5.16 &       0.00 &   11611.20 &          9 \\

         &      3 &       6.70 &       2.05 &       3.85 &       0.01 &    7147.18 &       5395 &       3.39 &       7.57 &       4.65 &       0.00 &    8616.02 &          1 &       2.03 &       9.83 &       5.13 &       0.00 &    9516.57 &          9 \\

\hline
\end{tabular}
\end{table}

\end{landscape}

We derive several key observations from our computational experiments. First, we compare the solution methods with respect to KPIs. 
In terms of ongoing corrective maintenance actions, {the \texttt{Stochastic Model} yields 1.79   jobs whereas \texttt{Strict Deterministic Model} and \texttt{Relaxed Deterministic Model} produce 5.69 and 6.59   jobs on average, respectively.} 
This suggests that  the \texttt{Stochastic Model} outperforms its deterministic counterparts {significantly} by proactively scheduling maintenance actions before failures occur. 
{Among deterministic models, the \texttt{Strict Deterministic Model}  model surpasses \texttt{Relaxed Deterministic Model} slightly but overall, both models fail to anticipate failures in a timely manner.}
As mentioned, the better performance of the \texttt{Stochastic Model} in terms of ongoing corrective jobs comes with an increase in  ongoing preventive maintenance jobs, {on average 5.38 jobs},  followed by the \texttt{Strict Deterministic} and \texttt{Relaxed Deterministic} models {by 1.33 and 0.77 jobs, respectively}.

The increase in ongoing preventive maintenance jobs not only decreases the corrective jobs but also reduces SLA violations. This is due to the shorter maintenance durations required for preventive jobs compared to corrective ones.  {In particular, the \texttt{Stochastic Model} results in 1.27 SLA  violations on average. On the other hand, \texttt{Strict Deterministic Model}  and \texttt{Relaxed Deterministic Model}  yield 1.61 and 1.78 SLA    violations on average, respectively. Therefore, the \texttt{Stochastic Model} is more advantageous than the deterministic models concerning this KPI.} We note that the average number of track violations is negligible for all the methods considered, and this KPI is omitted from the discussion below. 



The total cost analysis further supports our previous findings. Prioritizing the prevention of unexpected failures is rewarded most effectively by the \texttt{Stochastic Model} with a cost of 4838.12	 on average while average cost figures of the \texttt{Strict Deterministic Model} and the  \texttt{Relaxed Deterministic Model} are    5879.83  and   6279.73, respectively. This stark difference underscores the value of our approach. {Finally, total CPU time analysis suggest that \texttt{Stochastic Model} requires much more CPU time than its deterministic counterparts owing to the number of replications and large in-sample size.}

Next, we investigate the impact of {SLA} on the system. We note that the total weekly demand is 1809 for Line M4 whereas it is 1773 for Line M1B, indicating that Line M4 requires 2.03\% more railcars than Line M1B. This situation has some interesting consequences.
%
%
For example, the number of ongoing preventive and corrective maintenance jobs for Line M4 are larger than those of Line M1B   in all cases considered. 
 In addition, Line M4 also has higher SLA violations,  costs, and longer {CPU time} 
 compared to Line M1B. We note that higher values reported for Line M4 are possibly linked to higher total railcar demand {and fleet size as a result of it}.

Regarding fleet age, we analyze KPI trends across young, mixed, and old fleet compositions. We observe that an increase in fleet age   results in a higher number of ongoing maintenance jobs (both preventive and corrective) along with higher SLA violations and  total cost metrics in the majority of the settings considered.
Finally, we analyze the impact of gradually shifting the optimization priority from SLA violations to corrective maintenance by decreasing the $C_s / C_s$ ratio, in other words,  looking at Case 1, Case 2 and Case 3, respectively. {As the emphasis on corrective jobs increases, their number decreases in all settings while SLA violations and the number of preventive have a non-decreasing  trend. We observe lower overall costs {and CPU time in most cases} as we further prioritize corrective maintenance.} 

\begin{table}[H]
\centering
\caption{Extreme KPI metrics and the corresponding configurations. 
}
\label{tab:kpi-extremes}
\footnotesize
\begin{tabular}{|c|c|c|c|c|c|}
\hline
 \textbf{KPI}  &\textbf{Max/Min} &\textbf{Line} & \textbf{Fleet} & \textbf{Case}& \textbf{Value} \\
\hline
  \multirow{2}{*}{\#Prev}
  &Max &M4   & Old&   3   & 6.70 \\
  &Min &M1B   & Old &   3&4.65    \\
  \hline
  \multirow{2}{*}{\#Cor}        
  &Max &M4   & Mixed &   1 & 3.23    \\
  &Min &M1B  &  Young &   3 & 0.47    \\
            \hline
  \multirow{2}{*}{SLA\_v}     
  &Max &M4   & Old&   3   & 3.85   \\
  &Min    & M1B & Young & 1, 2 & 0.02     \\
             \hline
  \multirow{2}{*}{Cost}             
  &Max &M4  & Old &  1 & 10737.90  \\
  &Min &M1B   & Young & 3   & 2544.51  \\
          \hline
  \multirow{2}{*}{Time(s)}         
  &Max &M4   & Mixed  & 1 & 54030    \\
  &Min&M1B  & Young &   3  & 1350      \\
\hline
\end{tabular}
\end{table}

Table \ref{tab:kpi-extremes} summarizes the system configurations at which KPIs attain their extreme values for the \texttt{Stochastic Model}. {A young fleet for Line M1B  yields the smallest values for ongoing corrective maintenance jobs, SLA  violations, total cost and time. 
In contrast, an old fleet for Line M4 attains the highest level of SLA   violations in addition to the largest overall cost and the number of ongoing preventive maintenance jobs. Finally, a mixed fleet for Line M4  has the largest values for time and number of ongoing corrective jobs, which may stem from the fact that the \texttt{Stochastic Model} is trying to find a balance between preventive and corrective maintenance.}

\subsection{Confidence Interval Analysis}
\label{sec:compCI}

The out-of-sample analysis from Subsection~\ref{sec:compOutOfSample} only considers the point estimate of the key performance indicators. However, the SAA methodology used in the \texttt{Stochastic Model} can also provide confidence intervals for the statistical lower bound (LB) and upper bound (UB) of these indicators. Moreover, using the lower-end of LB and the upper-end of UB, we can compute a percentage gap metric, which statistically estimates a pessimistic bound on the relative optimality gap. On the other hand, deterministic models do not provide any statistical bounds on LB, and, consequently, any bound on the optimality gap.

We report the CI of total cost obtained from our experiments in Tables~\ref{tab:SAA-lineM1b} and \ref{tab:SAA-lineM4} for Lines M1B and M4, respectively. For deterministic models, we only report the CI for UB whereas for the \texttt{Stochastic Model}, we report the CIs for both LB and UB together with the percentage gap. It is evident that the \texttt{Stochastic Model} provides the best upper bounds, followed by the  \texttt{Strict Deterministic Model}.
In all cases and settings, we observe that the CIs are relatively small under all methods, confirming our conclusions based on point estimates from Subsection~\ref{sec:compOutOfSample}. In addition, the percentage gap values are less than {2.5\%} for the majority of case-setting combinations for the \texttt{Stochastic Model} and {the exceptions are no more than 3.2\%}. This suggests that our sample size selections are also quite reasonable for obtaining reliable results.

\begin{table}[H]
\footnotesize
\centering
\caption{Confidence intervals of the total cost for Line M1B.}
\label{tab:SAA-lineM1b}
\begin{tabular}{|c|c|c|c|c|c|c|}
\hline
\multicolumn{2}{|c|}{\textbf{Setting}} & \multicolumn{3}{c|}{\texttt{Stochastic Model}} & \texttt{Strict Deterministic} & \texttt{Relaxed Deterministic} \\
\hline
\textbf{Fleet} & \textbf{Case} & \textbf{CI for LB} & \textbf{CI for UB} & \textbf{\%Gap} & \textbf{CI for UB} & \textbf{CI for UB} \\
\hline
      \multirow{3}{*}{Young}  &      1 & [2503.48, 2511.93] & [2552.52, 2565.57] &       2.42 & [2936.29, 2964.74] & [2936.29, 2964.74] \\

       &      2 & [2503.06, 2510.53] & [2544.15, 2555.36] &       2.05 & [2913.01, 2937.54] & [2913.01, 2937.54] \\

       &      3 & [2501.04, 2508.79] & [2539.15, 2549.87] &       1.92 & [2896.57, 2918.64] & [2896.57, 2918.64] \\
\hline

      \multirow{3}{*}{Mixed}  &      1 & [3033.62, 3077.83] & [3100.96, 3133.98] &       3.20 & [4129.97, 4179.03] & [4504.49, 4557.66] \\

       &      2 & [2930.14, 2964.39] & [2997.31, 3021.66] &       3.03 & [3913.68, 3950.12] & [4215.52, 4254.39] \\

       &      3 & [2843.26, 2868.90] & [2887.79, 2907.01] &       2.19 & [3762.10, 3790.06] & [4013.08, 4042.26] \\
\hline

        \multirow{3}{*}{Old}  &      1 & [7322.37, 7371.95] & [7404.04, 7454.21] &       1.77 & [9039.99, 9115.08] & [9675.30, 9752.10] \\

         &      2 & [6010.15, 6044.98] & [6062.56, 6098.67] &       1.45 & [7356.36, 7412.59] & [7877.97, 7934.26] \\

         &      3 & [5081.07, 5108.68] & [5119.34, 5146.08] &       1.26 & [6179.43, 6221.98] & [6614.71, 6657.32] \\

\hline

\end{tabular}
\end{table}

\begin{table}[H]
\footnotesize
\centering
\caption{Confidence intervals of the total cost for Line M4.}
\label{tab:SAA-lineM4}
\begin{tabular}{|c|c|c|c|c|c|c|}
\hline
\multicolumn{2}{|c|}{\textbf{Setting}} & \multicolumn{3}{c|}{\texttt{Stochastic Model}} & \texttt{Strict Deterministic} & \texttt{Relaxed Deterministic} \\
\hline
\textbf{Fleet} & \textbf{Case} & \textbf{CI for LB} & \textbf{CI for UB} & \textbf{\%Gap} & \textbf{CI for UB} & \textbf{CI for UB} \\
\hline

      \multirow{3}{*}{Young} &      1 & [3337.60, 3353.96] & [3392.20, 3409.44] &       2.11 & [3798.33, 3834.66] & [3798.33, 3834.66] \\

       &      2 & [3331.38, 3349.77] & [3374.17, 3389.22] &       1.71 & [3744.71, 3775.00] & [3744.71, 3775.00] \\

       &      3 & [3325.79, 3343.07] & [3369.65, 3383.07] &       1.69 & [3707.08, 3733.34] & [3707.08, 3733.34] \\
\hline

      \multirow{3}{*}{Mixed} &      1 & [5447.06, 5535.62] & [5564.14, 5623.49] &       3.14 & [6492.89, 6554.74] & [7381.17, 7448.45] \\

       &      2 & [4865.79, 4930.59] & [4950.45, 4991.07] &       2.51 & [5798.14, 5843.58] & [6508.54, 6557.81] \\

       &      3 & [4427.98, 4476.39] & [4497.60, 4527.42] &       2.20 & [5312.64, 5346.96] & [5897.52, 5934.53] \\
\hline

        \multirow{3}{*}{Old} &      1 & [10587.30, 10644.70] & [10709.50, 10766.30] &       1.66 & [13967.30, 14051.90] & [14870.90, 14956.80] \\

         &      2 & [8531.70, 8569.99] & [8624.73, 8665.18] &       1.54 & [10900.00, 10963.40] & [11579.40, 11643.00] \\

         &      3 & [7066.93, 7096.21] & [7132.39, 7161.98] &       1.33 & [8592.03, 8640.01] & [9492.60, 9540.53] \\

\hline
\end{tabular}
\end{table}

\subsection{Detailed Solution Comparison}
\label{sec:detailedSolComp}

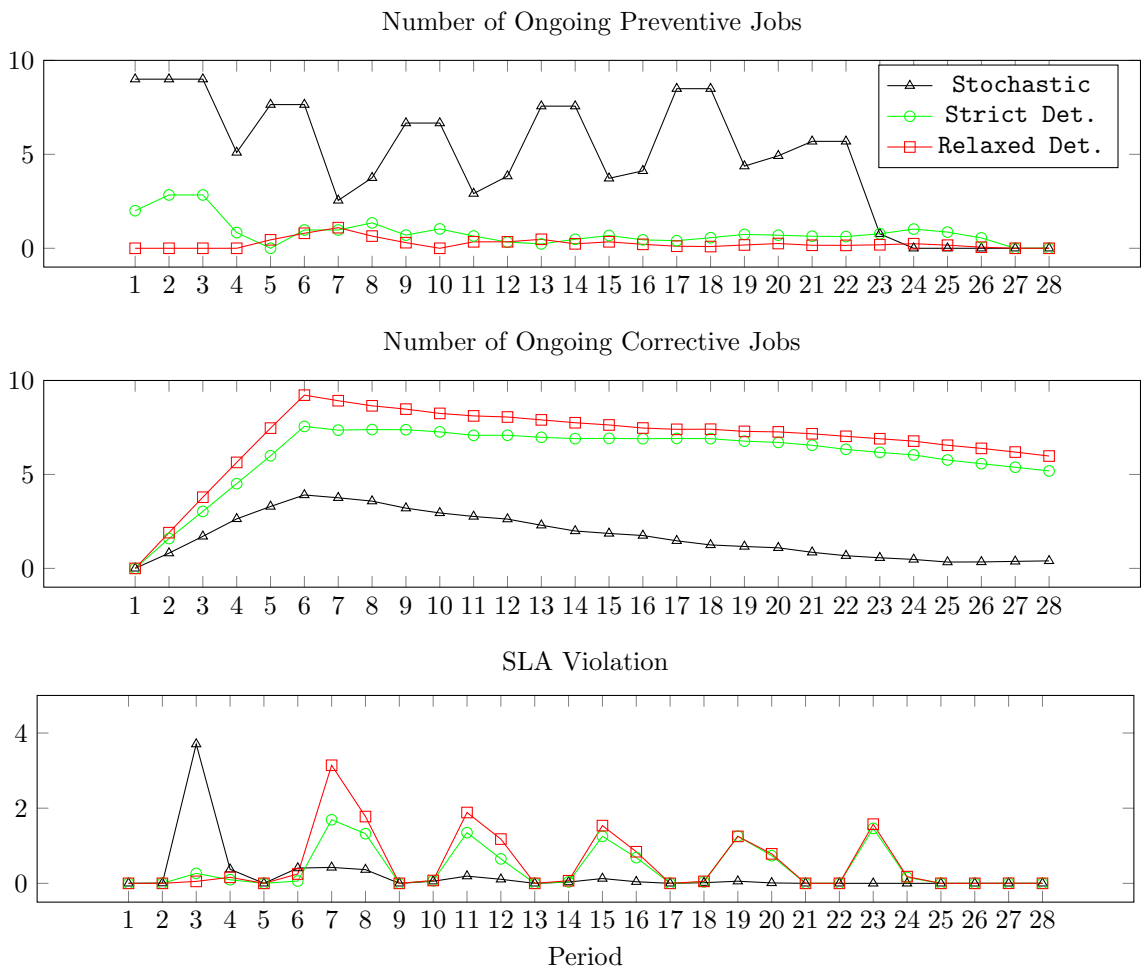
\begin{figure}[h]
\caption{The progression of three KPIs over the planning horizon for a mixed fleet in Line M1B under Case 1.
}
\label{fig:detailedComparison}
\centering
\pgfplotstableread[col sep=&, header=true]{
time	&	stoch	&	strict	&	relax
1	&	9	&	2	&	0
2	&	9	&	2.838	&	0
3	&	9	&	2.838	&	0
4	&	5.099	&	0.838	&	0
5	&	7.647	&	0	&	0.449
6	&	7.647	&	0.974	&	0.799
7	&	2.548	&	0.974	&	1.098
8	&	3.752	&	1.352	&	0.649
9	&	6.664	&	0.7	&	0.299
10	&	6.664	&	1.031	&	0
11	&	2.912	&	0.653	&	0.344
12	&	3.841	&	0.331	&	0.344
13	&	7.564	&	0.232	&	0.478
14	&	7.564	&	0.478	&	0.235
15	&	3.723	&	0.68	&	0.344
16	&	4.121	&	0.448	&	0.21
17	&	8.493	&	0.404	&	0.109
18	&	8.493	&	0.562	&	0.092
19	&	4.372	&	0.739	&	0.176
20	&	4.921	&	0.696	&	0.253
21	&	5.69	&	0.641	&	0.161
22	&	5.69	&	0.623	&	0.155
23	&	0.769	&	0.773	&	0.186
24	&	0	&	1.021	&	0.244
25	&	0	&	0.862	&	0.166
26	&	0	&	0.553	&	0.058
27	&	0	&	0	&	0
28	&	0	&	0	&	0
}\prevdata
     \begin{tikzpicture}
\begin{axis}[
    title={Number of Ongoing Preventive Jobs},
   xtick=data,
   xticklabels from table={\prevdata}{time},       
   ymax=10,  
     y=2.5mm, x=4.5mm,
           legend entries={\texttt{Stochastic}, \texttt{Strict Det.}, \texttt{Relaxed Det.} },
                    legend style={
        legend columns=1},
  ]
\addplot[mark=triangle] table [ y=stoch, x expr=\coordindex,] {\prevdata}; 
\addplot[ mark=o, color=green] table [y=strict, x expr=\coordindex,dashed] {\prevdata}; 
\addplot[ mark=square, color=red] table [y=relax, x expr=\coordindex] {\prevdata}; 
\end{axis}
    \end{tikzpicture}
\pgfplotstableread[col sep=&, header=true]{
time	&	stoch	&	strict	&	relax
1	&	0	&	0	&	0
2	&	0.812	&	1.588	&	1.901
3	&	1.705	&	3.031	&	3.788
4	&	2.634	&	4.509	&	5.64
5	&	3.295	&	5.993	&	7.465
6	&	3.907	&	7.555	&	9.22
7	&	3.758	&	7.358	&	8.923
8	&	3.579	&	7.388	&	8.652
9	&	3.205	&	7.381	&	8.476
10	&	2.951	&	7.264	&	8.248
11	&	2.763	&	7.08	&	8.115
12	&	2.623	&	7.085	&	8.056
13	&	2.296	&	6.976	&	7.903
14	&	1.99	&	6.906	&	7.751
15	&	1.859	&	6.917	&	7.633
16	&	1.749	&	6.902	&	7.469
17	&	1.468	&	6.914	&	7.4
18	&	1.248	&	6.904	&	7.407
19	&	1.166	&	6.771	&	7.291
20	&	1.095	&	6.703	&	7.264
21	&	0.856	&	6.552	&	7.163
22	&	0.674	&	6.333	&	7.029
23	&	0.565	&	6.172	&	6.899
24	&	0.476	&	6.038	&	6.773
25	&	0.337	&	5.767	&	6.553
26	&	0.343	&	5.572	&	6.387
27	&	0.373	&	5.383	&	6.191
28	&	0.4	&	5.185	&	5.979
}\cordata
     \begin{tikzpicture}
\begin{axis}[
    title={Number of Ongoing Corrective Jobs},
   xtick=data,
   xticklabels from table={\cordata}{time},      
   ymax=10,
     y=2.5mm, x=4.5mm
  ]
\addplot[mark=triangle] table [ y=stoch, x expr=\coordindex,] {\cordata}; 
\addplot[ mark=o, color=green] table [y=strict, x expr=\coordindex,dashed] {\cordata}; 
\addplot[ mark=square, color=red] table [y=relax, x expr=\coordindex] {\cordata}; 
\end{axis}
    \end{tikzpicture}
\pgfplotstableread[col sep=&, header=true]{
time	&	stoch	&	strict	&	relax
1	&	0	&	0	&	0
2	&	0.006	&	0	&	0
3	&	3.705	&	0.265	&	0.052
4	&	0.375	&	0.096	&	0.16
5	&	0	&	0	&	0
6	&	0.406	&	0.061	&	0.248
7	&	0.424	&	1.692	&	3.142
8	&	0.365	&	1.322	&	1.777
9	&	0	&	0	&	0
10	&	0.056	&	0.069	&	0.076
11	&	0.192	&	1.349	&	1.885
12	&	0.108	&	0.651	&	1.178
13	&	0	&	0	&	0
14	&	0.034	&	0.038	&	0.069
15	&	0.128	&	1.251	&	1.536
16	&	0.043	&	0.684	&	0.841
17	&	0	&	0	&	0
18	&	0.017	&	0.039	&	0.052
19	&	0.057	&	1.256	&	1.249
20	&	0.011	&	0.736	&	0.785
21	&	0	&	0	&	0
22	&	0	&	0	&	0
23	&	0	&	1.463	&	1.575
24	&	0	&	0.157	&	0.179
25	&	0	&	0	&	0
26	&	0	&	0	&	0
27	&	0	&	0.001	&	0.002
28	&	0	&	0	&	0
}\sladata
     \begin{tikzpicture}
\begin{axis}[
    title={SLA Violation},
   xtick=data,
   xticklabels from table={\sladata}{time},      
   ymax=5,   
			    xlabel={Period},
     y=5mm, x=4.5mm,
  ]
\addplot[mark=triangle] table [ y=stoch, x expr=\coordindex,] {\sladata}; 
\addplot[ mark=o, color=green] table [y=strict, x expr=\coordindex,dashed] {\sladata}; 
\addplot[ mark=square, color=red] table [y=relax, x expr=\coordindex] {\sladata}; 
\end{axis}
    \end{tikzpicture}
\end{figure}

This section aims to give further insight as to how the \texttt{Stochastic Model} outperforms the \texttt{Strict Deterministic Model}  and the  \texttt{Relaxed Deterministic Model}. For this illustration, we have chosen the setting in which the relative cost improvement of the \texttt{Stochastic Model} is the highest: A mixed fleet in Line M1B under Case 1. We report the progression of three KPIs over the planning horizon in Figure~\ref{fig:detailedComparison} for the three methods.

We observe that the \texttt{Stochastic Model}  schedules more preventive jobs than its deterministic competitors, especially in the beginning of the planning horizon. 
These predictive maintenance decisions pay off in terms of reducing the number of corrective jobs. As it can be seen from the second subfigure in Figure~\ref{fig:detailedComparison}, both the \texttt{Strict Deterministic Model} and the \texttt{Relaxed Deterministic Model} suffer from  corrective jobs through the planning horizon whereas the \texttt{Stochastic Model} has significantly smaller number of failures. Consequently, the SLA violation is mostly kept around zero for the \texttt{Stochastic Model}, except for period 3, at which point the model has already scheduled a large number of preventive jobs. 

Overall, we can clearly see that the predictive maintenance decisions taken in the beginning of the planning horizon by the \texttt{Stochastic Model} help keep both railcar failures and SLA violations much smaller compared to \texttt{Strict Deterministic} and \texttt{Relaxed Deterministic} models. Consequently, the overall cost of the \texttt{Stochastic Model}  is significantly lower.

\section{Conclusions} \label{s:conc}

In this study, we address predictive railcar-maintenance scheduling by jointly considering preventive and corrective maintenance costs, SLA compliance, and track-capacity constraints. We propose a two-stage stochastic programming framework solved with the SAA method and benchmark it against two deterministic models that differ only in their preventive-maintenance interval policies.

The out-of-sample evaluation shows that the \texttt{Stochastic Model} anticipates clustered failures by withdrawing railcars earlier for preventive maintenance, thereby reducing corrective jobs and SLA violations. This translates into a 17.72\% cost improvement, on average, over the best deterministic model, albeit at the expense of substantially higher computational effort. Among deterministic approaches, the strict policy proves more effective than the relaxed one, achieving a 6.37\% cost reduction. Confidence interval analysis further supports the robustness of our results, with the pessimistic cost gap never exceeding 3.2\%.

System parameters also shape KPI outcomes. Line M4, having 2.03\% more total demand and operating with 29.82\% more fleet size, naturally requires more maintenance, incur more SLA violations. Besides, older fleets require more maintenance and incur more SLA violations, while mixed fleets lead to the largest CPU times and number of ongoing corrective jobs, reflecting the additional complexity of balancing preventive and corrective jobs. 

There are promising future research directions. In terms of application, a model handling the railcars causing track capacity violation as queue and choosing longer planning horizon to capture its carryover effect to the next scheduling periods would be more realistic. Secondly, our experiments suggest that \%3.20 gap is noticed in SAA method at worst, so to further diminish it, we need to increase in-sample size at the expense of highly CPU time observed at our preliminary computations. To solve this issue, Benders Decomposition method can be utilized.

\section*{Data Statement}

The codes are available at https://github.com/muratelhuseyni/Stochastic-Railcar-Maintenance

\section*{Statements and Declarations}

\subsection*{Competing interests}
The authors declare that they have no competing interests.







\subsection*{Availability of data and materials}
The data associated with the manuscript is available upon reasonable request. 
The code is available at \url{https://github.com/muratelhuseyni/Stochastic-Railcar-Maintenance}.





\bibliographystyle{plain}
\bibliography{references}

\end{document}